\documentclass[12pt]{article}
\usepackage{latexsym}
\usepackage{amssymb}
 
\newcommand{\proof}{{\noindent \bf Proof. }}
\newtheorem{thm}{Theorem}
\newtheorem{defi}{Definition}
\newtheorem{lem}{Lemma}

\newtheorem{prop}{Proposition}

\newcommand{\N}{{\mathbb N}}

\newcommand{\C}{{\cal C}}

\newcommand{\F}{{\cal F}}

\newcommand{\K}{{\cal K}}
\newcommand{\U}{{\cal L}}

\date{}

\begin{document}
\begin{titlepage}
\title{\bf ON THE PERMUTATION CAPACITY OF DIGRAPHS}
{\author{{\bf  G\'erard Cohen}
\\{\tt cohen@enst.fr}
\\''ENST
\\ FRANCE
\and{\bf Emanuela Fachini}
\\{\tt fachini@di.uniroma1.it}
\\''La Sapienza'' University of Rome
\\ ITALY
\and{\bf J\'anos K\"orner}
\thanks{Department of Computer Science, University of Rome, La Sapienza, 
via Salaria 113, 00198 Rome, ITALY}
\\{\tt korner@di.uniroma1.it}
\\''La Sapienza'' University of Rome
\\ ITALY}} 

\maketitle
\begin{abstract}

We extend several results of the third author and C. Malvenuto on graph--different permutations to the case of directed graphs and introduce new open problems. Permutation capacity is a natural extension of Sperner capacity from finite directed graphs to infinite digraphs. Our subject is combinatorial in nature, but can be equally regarded as zero--error information theory.

\end{abstract}
\end{titlepage}

\section{A puzzle}

How many permutations of the first $n$ natural numbers can we find such that any two of them place 
two consecutive natural numbers somewhere in the same position? This problem was introduced in
\cite{KM} where the authors conjectured that the maximum number of such permutations is exactly the middle binomial coefficient 
${n \choose {\lceil {n \over 2} \rceil}}.$ The conjecture has been neither disproved nor confirmed so far. This combinatorial  puzzle is closely connected to Shannon's graph capacity concept and the rich and beautiful mathematics around it. In this paper we begin to explore its generalizations to directed graphs. 

\section{Introduction for the information theorist}

In \cite{CK} the zero--error capacity of a discrete memoryless stationary channel was generalized by restricting the input sequences to "mimick" a fixed distribution. The resulting concept of {\it zero--error capacity within a given type} 
allows to apply the method of types \cite{CsK} to several important problems in extremal combinatorics and the theory of combinatorial search. All such problems can be regarded as part of zero--error information theory. In \cite{KS} and the follow--up paper \cite{GKV} Shannon's graph capacity problem \cite{Sh} was generalized to directed graphs. The resulting concept of 
Sperner capacity and the corresponding capacity within a given type gave the key to solve an intriguing problem of 
R\'enyi in combinatorial search in \cite{GGKV}. 
This direction of research started in \cite{CKS} where the problem of 
zero--error capacity of the compound channel with uninformed encoder (but not decoder) was introduced as a unifying model for several problems in extremal combinatorics. Actually, in 
\cite{GGKV} a formally information--theoretic analogue of this problem was solved for families of directed graphs. Much to our surprise this generalization gave 
Nayak and Rose \cite{NR} the correct mathematical formulation and the technique of solution for the important problem of the zero--error capacity of the compound channel in case of uninformed encoder and decoder.

In this paper we introduce and study corresponding problems for infinite alphabet channels. 
Formally, our problems are concerned with Sperner--type capacity of infinite directed graphs. To the analogy of the case of capacity within a given type we restrict the input sequences to be permutations. In case of undirected graphs these problems have been 
introduced in \cite{KM} and further studied in the papers  \cite{KMS} and \cite{KSS}.

\section{Digraph--different permutations}

Let $\N$ denote the set of natural numbers and let $G$ be an arbitrary directed graph with vertex set $\N$. We will say that two permutations, $\pi$ and $\rho$ of the first $n$ natural numbers are 
$G$--{\em different} if there is an $i \in [n]$ such that the ordered couple 
of its images under these two permutations satisfies $(\pi(i), \rho(i)) \in E(G).$
We write $N(G, n)$ for the largest cardinality of a set of pairwise $G$--different permutations of $[n]$. 
It is easy to see (cf. \cite{KMS} for the analogous statement for undirected graphs) that
$N(G, n)$ is supermultiplicative in $n$ for any fixed $G$, and thus 
the expression $\frac{1}{n} \log N(G,n)$ has a limit, even though its value is not always finite. 
We denote by $R(G)$ the always existing limit
$$R(G)=\lim_{n \rightarrow \infty} \frac{1}{n} \log N(G,n)$$
and call it the {\it permutation capacity} of $G.$
All this is in perfect analogy with the undirected case. In fact, our present concepts 
are formal generalizations of what has been introduced in previous work, mainly 
in \cite{KMS} and \cite{KSS}, for undirected graphs. As it is well--known and easily seen, any undirected graph $F$
can be formally identified with the {\it symmetrically directed graph} $G_F$ on the same vertex set for which $(a,b)\in E(G_F)$ if and only if $\{a,b\}\in E(F).$
Here and in the sequel all the logarithms are to the base 2.
We observe that
\begin{lem}\label{lem:chro}
$$R(G)\leq \log \chi(G)$$
where $\chi(G)$ is the chromatic number of the undirected graph underlying $G$.
\end{lem}

\proof

Clearly, if for two digraphs $G_1$ and $G_2$ on the same vertex set $\N$ we have
$E(G_1)\subseteq E(G_2)$ then $R(G_1)\leq R(G_2).$ 
This implies that $R(G)\leq R(G')$ where $G'$ is the digraph having two, oppositely 
oriented edges for every edge of $G.$ We observe that, by definition,  $R(G')$ is equal to the permutation capacity of the undirected graph underlying  $G'.$ The rest follows by applying to this undirected graph (the fairly obvious) Proposition 4 from \cite{KSS}.
\hfill$\Box$
 
We will concentrate on the various digraphs whose underlying undirected graph is the infinite path $L.$ In this graph  
every pair of consecutive integers is adjacent. By the foregoing, all the directed graphs 
defined on $L$ have capacity at most 1. 
(The classical Shannon capacity of an infinite graph is not usually defined even though it makes perfect sense and just as in the finite case is trivial to determine if the graph in question has, as in our present case, coinciding clique number and chromatic number.)
Obviously, the maximum of these capacities 
is achieved by the {\it symmetric path} $L_{sym}$ whose edge set contains two oppositely oriented edges for every edge of $L.$ It was shown in \cite{KMS} that the 
permutation capacity of this symmetric graph is at least $\log\sqrt[4]{10}.$ (This bound was further strengthened by Brightwell and Fairthorne \cite{BF}).
No non--trivial upper bound for $R(L_{sym})$ is known. In \cite{KM} the trivial upper bound 1 was conjectured to be tight. 

New interesting questions arise if we concentrate on {\it oriented path}s in which every couple of vertices gives rise to at most one oriented edge. Let $\U$ denote the infinite family of all the oriented paths so obtained.
Our main result is the following inequality:
$$\frac{1}{2}< \inf_{G \in \U} R(G)$$

We are especially interested in 
two simple oriented paths defined for $L.$ We will denote by $L_c$ the digraph in which every edge of $L$ is oriented from its smaller vertex to the larger one and call it the 
{\it thrupath}. Similarly, we will denote by $L_a$ the {\it alternating path} in which every 
edge is oriented from its odd vertex towards its even vertex. 
We will prove in this paper that  $N(L_a, n)$ growths with $n$ as a standard Fibonacci sequence while $N(L_c, n)$ exhibits the same growth rate in an asymptotic sense.  

\section{Diverse capacities} 

The study of the Sperner capacity of digraphs was prompted by an interest in solving a famous problem of R\'enyi about the largest cardinality of a family of pairwise qualitatively independent $k$--partitions of an $n$--set.  For this purpose, Simonyi and the third author initially proposed to study a concept of  
capacity for undirected graphs \cite{KS} that later turned out to be reducible to the subsequent more subtle concept of Sperner capacity, \cite{GGKV}, \cite{Cfive}. This concept seems however less tractable and more intriguing in case of the permutation capacity of infinite graphs. 

\begin{defi}\label{defi:ks}
Let $G$ be an undirected graph with vertex set $\N$. We will say that the permutations $\pi$ and 
$\rho$ of $[n]$ are {\it robustly $G$--different} if there are two elements $i\in [n]$ and $j \in [n]$ 
such that $(\pi(i), \rho(i))=(\rho(j), \pi(j))$ and $\{\pi(i), \rho(i)\} \in E(G).$

Let $NN(G, n)$ be the maximum cardinality of a set of pairwise robustly $G$--different permutations of 
$[n].$ We call the limit
$$RR(G)=\lim_{n \rightarrow \infty} \frac{1}{n} \log NN(G,n)$$
the robust permutation capacity of $G.$
\end{defi}

It is interesting to relate this quantity to the infimum $R_{min}(L)$ of the 
permutation capacities of all the digraphs obtained from $L.$ One might conjecture that, as in the case of the finite family of oriented versions of a finite simple graph, this 
infimum equals the robust permutation capacity of $L.$ 
It is very easy to see (but we will return to this) that

$$RR(L)\geq \frac{1}{2}$$

We conjecture this lower bound to be tight. 
If true, this conjecture, when compared to our Theorem \ref{thm:one} below,  would imply that the main theorem of \cite{GGKV} does not generalize to infinite families of graphs. 

We consider analogous problems for other digraphs on $\N$ including infinite tournaments. In case the corresponding capacity is infinite, we ask more refined 
questions about the rate of asymptotic growth of $N(G,n).$ 

\section{Constructions}

To begin, we will study the particularly simple case of the alternating path.
Consider the standard Fibonacci sequence $f(1)=1,$ $f(2)=2$, $f(n)=f(n-1)+f(n-2).$ We claim

\begin{prop}\label{prop:fib}
$$N(L_a, n)\geq f(n)$$
\end{prop}

\proof
Let $\sigma_i$ be the transposition $(i, i+1).$ Then for every set $F\subseteq [n-1]$ of natural numbers 
not containing consecutive elements consider the following product of 
transpositions
$$\sigma_F=:\prod_{i \in F} \sigma_i$$  (As all these transpositions commute, the product is well--defined.) We claim that as 
$F$ runs over all the subsets of $[n-1]$ without consecutive elements, the permutations so obtained give the desired construction.

In what follows we will call a subset of $\N$ {\it loose} if it does not contain two consecutive integers.
Let $E$ and $F$ be two different loose subsets of $[n-1]$ and let $i$ be the smallest integer in  $[n]$
contained in exactly one of the two. Without loss of generality, suppose that $i \in E$, $i \not \in F.$
Then $i-1$ is not contained in any of the two sets. 
This means that 
$\sigma_F(i)=i$ while $\sigma_E(i)=i+1.$ 
Suppose for a moment that $i$ is odd.  This implies that
$$(\sigma_F(i), \sigma_E(i)) \in E(L_a).$$
In particular, note that in the two different linear orderings of $[n]$ the different permutations 
$\sigma_F$ and $\sigma_E$ define, the subsets of $[n]$ determined by the positions of the even (odd) integers in $[n]$ are different. 
This means that necessarily, for some other value $j \in [n]$ the previous relation has to be "reversed" at least in the sense that $\sigma_F(j)$ is even and $\sigma_E(j)$ is odd. On the other hand, since both of our permutations have the property that 
$$|\sigma(i)-i|\leq 1$$
for every $i \in [n],$ we see that
$$|\sigma_F(j)-\sigma_E(j)|\leq 2.$$ 
Since we know that the two images of $j$ have different parity, we must conclude that 
$$|\sigma_F(j)-\sigma_E(j)|=1$$ 
and thus 
$$(\sigma_E(j), \sigma_F(j)) \in E(L_a).$$
This concludes the proof that our permutations differ in the required way under the hypothesis that the first position in which $E$ and $F$ differ is 
odd. Clearly, if that position is even, the same proof can be applied exchanging odd and even in the above.

Finally, note that the number of loose subsets of $[n-1]$ (which is also the cardinality of our construction),  is exactly $f(n).$ 
\hfill$\Box$

We conjecture that the alternating path yields the maximum of the permutation capacities of the 
digraphs in $\U.$
For the thrupath we have a similar construction, yielding a slightly weaker lower bound.
\begin{prop}\label{prop:thru}

Let $L_c$ denote the thrupath. Then 
$$N(L_c,n)\geq \frac{f(n)}{n^3}$$
whence
$$R(L_c)\geq \log \frac{\sqrt{5}+1}{2}$$

\end{prop}

\proof

Let $L_c$ be the oriented graph obtained from the infinite path by orienting its edges from their smaller endpoint to the greater one. Fix an arbitrary natural number 
$n \in \N.$ Our construction has the same starting point as the one in the previous proposition. For any loose subset $F$ of $[n-1]$ we will consider the permutation $\sigma_F$ defined there. To each of these permutations we will associate the couple of values of two functionals we define next. 

Let $F:=\{f_1, f_2,\dots, f_k\}$ be an arbitrary loose subset of $[n-1]$ where $k=|F|.$ 
We define
$$\lambda(F):=|F| \quad
\hbox{and}
\quad \mu(F):=\sum_{i \in F}i.$$
These functionals have relatively small ranges as $F$ varies in $[n-1].$ In fact, $\lambda$ has at most $n$ different values, while the range of $\mu$ has less than $n^2$ elements.
Thus the couple of these functionals, considered as a single functional, has a range of less than $n^3$ elements. 
Now let us look at the partition of the family of loose subsets of $[n-1]$ generated by the full inverse images of the various values of the last functional.
By the foregoing, this partition has at most $n^3$ different classes. Hence there exists at least one class $\C$ 
whose cardinality gives
$$|\C|\geq \frac{f(n)}{n^3}.$$
Note that every set $F \in \C$ has the same cardinality. Let $k$ be this common cardinality. 

We claim that the permutations associated by our construction to the member sets of $\C$ are pairwise 
$L_c$-different. In order to see this, let us consider two different sets, $E$ and $F$ from $\C$ and write
$$E:=\{e_1, e_2, \dots, e_k\} \quad \hbox{and} \quad F:=\{f_1, f_2,\dots, f_k\}.$$ Let $j$ be the first 
index for which $e_j \not= f_j.$ Without loss of generality we can suppose that $e_j<f_j.$ This means for the corresponding permutations that
$\sigma_E(e_j)=e_j+1$ while $\sigma_F(e_j)=e_j,$ whence 
$$(\sigma_F(e_j), \sigma_E(e_j))\in E(L_c).$$
Now, by our construction
$$\sum_{i=1}^ke_i=\sum_{i=1}^kf_i.$$
This implies that there is an index $i$ for which $e_i>f_i.$ Let $\l$ denote the first such index. 
It is easy to see that we have 
$\sigma_F(f_i)=f_i+1,$ while at the same time $\sigma_E(f_i)=f_i.$ The latter is implied by the 
fact that $f_i-1 \not \in F$ (since $F$ is loose) and thus $f_i-1 \not \in E.$ This gives
$$(\sigma_E(f_i), \sigma_F(f_i))\in E(L_c).$$

\hfill$\Box$

To conclude this section, we will concentrate on the robust permutation capacity of $L.$

\begin{lem}\label{lem:rob}

For the infinite path $L$ we have
$$NN(L, n)\geq 2^{\lfloor {n \over 2} \rfloor}$$
implying
$$RR(L)\geq \frac{1}{2}$$
\end{lem}

\proof

Let $\F$ run through the subsets of the odd elements of $[n-1].$ Clearly 
$|\F|=2^{\lfloor {n \over 2} \rfloor}.$
On the other hand, consider any two sets from $\F.$ Let these be $E$ and $F$ and let 
$j$ be the first element contained in only one of them, say $j\in E-F.$ Then, as before, 
we have for the corresponding permutations 
$$\sigma_E(j)=j+1, \quad \sigma_E(j+1)=j \quad \hbox{while} 
\quad \sigma_F(j)=j, \quad \sigma_F(j+1)=j+1.$$  

\hfill$\Box$

It seems very interesting to establish whether this lower bound is tight, especially in the light 
of our main result, concerning the infimum of the permutation capacities of the digraphs 
from $\U.$

\begin{thm}\label{thm:one}
$$\inf_{G\in \U}R(G)\geq \frac{\log 3}{3}>\frac{1}{2}.$$
\end{thm}

\proof

We define the digraph  $B=\{V(G), E(G)\}$ as follows. The vertices of $B$ are 
$V(B):=\{a,b, c\}$ and the edges are 
$E(B):=\{(a,c), (c,a), (b,c), (c,b), (a,b)\}.$ In other words, $B$ is obtained from the symmmetric 
clique on 3 vertices by omitting a single edge. We denote by $\omega_s(B^n)$ 
the cardinality of the largest symmetric clique in the $n$'th conormal power of $B^n.$
It is well-known and easily seen (cf. \cite{GKV}) that 
\begin{equation}\label{eq:spe}
\omega_s(B^n)\geq \frac{3^n}{2n+1}
\end{equation}
for every $n \in \N.$

Let us now take a closer look at an arbitrary but fixed digraph $G \in \U.$
For this purpose, we consider, for every $n \in \N$ another digraph, $G_n$, whose vertex set is the set of permutations of $[n]$ and in which there is an edge  $(\rho, \sigma)$ pointing from the permutation
$\rho$ to $\sigma$ if there is a number $j \in [n]$ such that $(\rho(j), \sigma(j)) \in E(G).$ 
Clearly, we have 
$$N(G,n)=\omega_s(G_n).$$
We would like to show that
\begin{equation}\label{eq:main}
\limsup_{n \rightarrow \infty} \frac{1}{n}\log\omega_s(G_n)\geq \frac{\log 3}{3}
\end{equation}
This, however, will immediately follow from (\ref{eq:spe}) if we show that 
for every $n\in \N$ the graph $G_{3n}$ contains a subgraph isomorphic to $B^n.$
Here is how we achieve this objective. Consider first the vertices $\{1,2,3\}$ of our graph $G$ 
and identify the vertex set $V(B)$ of the digraph $B$ with permutations of 
$[3]=\{1,2,3\}$ as follows:
$$c:=(1,2,3),$$
and
 $$\quad a:=(1,3,2), \quad b:=(2,1,3) \quad \hbox{if both}  \quad 
(1,2) \in E(G)\quad \hbox{and} \quad (2,3) \in E(G)$$
while in the opposite case exchange the role of $a$ and $b$ by setting
$a:=(2,1,3)$ and $b:=(1,3,2).$ 
This procedure can be extended in the same way to any consecutive disjoint interval. 
More precisely, let us consider the vertex set $\N$ of $G$ as the disjoint union of the 
"intervals"
$\N=\bigcup_{t=0}^{\infty}\{3t+1, 3t+2, 3t+3\}$ and for every $t$ let us define the graph 
$B_t$ with vertex set 
$$c_t:=(3t+1, 3t+2, 3t+3),$$
and
$$ a_t:=(3t+1, 3t+3, 3t+2), \quad b_t:=(3t+2, 3t+1, 3t+3)$$ 
if both
$$(3t+1, 3t+2) \in E(G) \quad \hbox{and} \quad  (3t+2, 3t+3) \in E(G),$$
while in the remaining cases we define 
$$a_t:=(3t+2, 3t+1, 3t+3) \quad \hbox{and} \quad b_t:=(3t+1, 3t+3, 3t+2).$$
It is then clear that $B^n$ is isomorphic to the co--normal product of the graphs 
$B_t$ whence 
$$B^n\cong \bigotimes_{t=0}^{n-1}B_t\subseteq G_{3n}$$ 
and in virtue of (\ref{eq:spe}) this implies
$$\omega_s(G_{3n}) \geq \frac{3^n}{2n+1},$$
thereby completing our proof.

\hfill$\Box$

\par\noindent{\bf Remark}
\medskip

This statement is interesting because it cuts away the infimum of the permutation capacities for the family of graphs $\U$ from the conjectured robust capacity of the underlying path graph.Using our proof technique above, G\'abor Simonyi \cite{GS} 
has improved our lower bound to $\frac{\log 5}{4}.$ It is unlikely that even this
new bound be tight. Simonyi's improvement is based on the observation that
if one uses permutations of $[4]=\{1,2,3,4\}$ instead of those of $[3]=\{1,2,3\}$ then we can define, in the role of $B,$ a digraph whose vertices are the 5 permutations from our corresponding Fibonacci construction in the proof of Proposition \ref{prop:fib}. It is easy to see that this graph contains a transitively oriented clique containing all of its vertices. 
We leave the details of the proof to the reader. 
\section{Tournaments}
We would like to comment briefly on the apparently much simpler case of the family of all oriented 
versions of the infinite complete graph $K$ with vertex set $\N.$ Let us indeed  look at the family 
$\K$ of all the oriented graphs resulting from the different orientations of the edges of $K.$
In this case it is obvious that 
$N(F,n)$ takes its maximum value in $\K$ on the digraph $F \in \K$ only in case  the digraph induced 
by $F$ on $[n]$ 
is transitive. Surprisingly, it is considerably less obvious on which graph the minimum 
value is taken. One might conjecture, however, that if $n$ is a power of 3, $n=3^t,$ say, then this graph is the $t'$s lexicographic power of the so--called cyclic triangle. More precisely, we call cyclic triangle  the graph $C$ whose vertices $\{1,2,3\}$ are oriented cyclically, from 1 to 2, from 2 to 3 and then from 
3 to 1. Given a graph $C$, its $t$'th lexicographic power $C_{lexi}^{ t}$ is the graph with vertex set 
$[V(C)]^t$ in which $({\bf v,  w}) \in E(C_{lexi}^{ t})$ if, denoting by $j$ the first coordinate in which the 
strings ${\bf v}:=v_1, v_2, \dots, v_t$ and ${\bf w}:=w_1, w_2, \dots w_t$ differ, we have
$(v_j, w_j) \in E(C).$ It is well-known that this is the "most symmetric" tournament on $3^t$ vertices, 
having the largest automorphism group among all the digraphs with the same set of vertices; (this is the celebrated Alspach-Dixon theorem,  \cite{A} and \cite{D}.)
At the same time, this graph has the maximum number of induced 3--element  subgraphs which are cyclic triangles. However, as it is well--known,  
every tournament in which every vertex has the same out--degree and 
in--degree, also has this same number of cyclic triangle subgraphs.  

All the problems considered for simple graphs in \cite{KM}, \cite{KMS} and \cite{KSS} can be generalized 
to directed graphs and we plan to return to the corresponding new questions in subsequent research. In particular, it seems especially interesting to investigate the functional corresponding to $\kappa$ from 
\cite{KMS} in the case of finite digraphs. 
\section{Acknowledgement}

We are grateful to  Blerina Sinaimeri and G\'abor Simonyi for sharing their ideas with us.


\begin{thebibliography}{99}

\bibitem{A} B. Alspach, A combinatorial proof of a conjecture of Goldberg and Moon, 
{\it Canadian Math. Bull.}vol. 11 (1968), no. 5, pp. 655--661,

\bibitem{BF} G. Brightwell, M. Fairthorne, Permutation capacity of graphs, 
paper in preparation.

\bibitem{Cfive} R. Calderbank, P. Frankl, R. L. Graham, W. Li and L. Shepp, 
The Sperner capacity of the cyclic triangle for linear and nonlinear codes, 
{\it J. Algebraic Combin.}, {\bf 2}(1993), pp. 31--48,

\bibitem{CKS} G. Cohen, J. K\"orner, G. Simonyi, Zero-error capacities and very different sequences, 
in:{\it Sequences. Combinatorics, Security and Transimission},
Advanced International Workshop on Sequences, Positano, Italy, June 1988,
Springer, New York, 1990,  R. M. Capocelli, ed., 144--155,

\bibitem{CK} I. Csisz\'ar, J. K\"orner, On the capacity of the arbitrarily
  varying channel for maximum 
probability of error, {\it Zeitschrift f\"ur Wahrscheinlichkeitstheorie verw.
Geb}., {\bf 57} (1981), 87--101.

\bibitem{CsK} I. Csisz\'ar, J. K\"orner,
{\em Information theory: Coding theorems for discrete
memoryless systems}, Academic Press, New York, 1982 and Akad\'emiai
Kiad\'o, Budapest,

\bibitem{D} J. D. Dixon, The maximum order of the group of a tournament, 
{\it Canadian Math. Bull.}vol. 10(1967) no. 4, pp.503-505,

\bibitem{GKV} L. Gargano, J. K\"orner, U. Vaccaro, Qualitative independence and Sperner problems for directed graphs,
{\it J. Comb. Theory}, Ser. A, {\bf 61}(1992), 173--192,

\bibitem{GGKV} L. Gargano, J. K\"orner, U. Vaccaro,
Capacities: from information theory to extremal set theory,
{\it J. Comb. Theory} Ser. {\bf A}, {\bf 68}(1994), no. 2, pp. 296--315,

\bibitem{KM} J. K\"orner, C. Malvenuto, Pairwise colliding permutations and
the capacity of  infinite graphs, 
{\it SIAM J. Discrete Mathematics}, {\bf 20} (2006),
203--212. 

\bibitem{KMS}  J. K\"orner, C. Malvenuto, G. Simonyi, Graph--different
  permutations,  
{\it SIAM J. Discrete Mathematics},  2, vol. 22(2008), pp. 489--499,

\bibitem{KS} J. K\"orner, G. Simonyi, A Sperner-type theorem and qualitative independence, 
{\it J. Comb.Theory}, Ser. {\bf A}, 1({\bf 59})(1992), 90--103

\bibitem{KSS} J. K\"orner, G. Simonyi, B. Sinaimeri, 
On types of growth for graph--different permutations, 
{\it J. Comb. Theory Ser. A}, submitted,

\bibitem{NR} J. Nayak, K. Rose, Graph capacities and zero--error transmission over compound channels, 
{\it IEEE Trans. Inform. Theory}, {\bf 51}  no. 12(2005), 4374--4378,

\bibitem{Sh} C. E. Shannon, The zero--error capacity of a noisy channel, 
{\it IRE Trans. Inform. Theory}, {\bf 2} (1956), 8--19. 

\bibitem{GS} G. Simonyi, personal communication

\end{thebibliography}
\end{document}